\documentclass[a4paper,12pt]{article}
\usepackage{amscd}
\usepackage{amsmath,amsfonts,amssymb,amscd}
\usepackage{indentfirst,graphicx,epsfig}
\usepackage{graphicx,psfrag}
\input{epsf}

\newtheorem{thm}{Theorem}
\newtheorem{lem}{Lemma}
\newtheorem{con}{Conjecture}

\def\qed{\hfill \nopagebreak\rule{5pt}{8pt}}
\def\pf{\noindent {\it Proof.} }

\title{\bf Note on a relation between Randi\'c
index\\ and algebraic connectivity\footnote{Supported by NSFC and
``the Fundamental Research Funds for the Central Universities". } }
\author{
\small  Xueliang Li, Yongtang Shi\\
\small Center for Combinatorics and LPMC-TJKLC \\
\small Nankai University, Tianjin 300071, China \\
\small lxl@nankai.edu.cn,  shi@nankai.edu.cn
\date{}}
\begin{document}
\maketitle
\begin{abstract}
A conjecture of AutoGraphiX on the relation between the Randi\'c
index $R$ and the algebraic connectivity $a$ of a connected graph
$G$ is:
$$\frac R a\leq
 \left(\frac{n-3+2\sqrt{2}}{2}\right)/\left(2(1-\cos {\frac{\pi}{n}})\right)  $$
with equality if and only if $G$ is $P_n$, which was proposed by
Aouchiche and Hansen [M. Aouchiche and P. Hansen, A survey of
automated conjectures in spectral graph theory, {\it Linear Algebra
Appl.} {\bf 432}(2010), 2293--2322]. We prove that the conjecture
holds for all trees and all connected graphs with edge connectivity
$\kappa'(G)\geq 2$, and if $\kappa'(G)=1$,  the conjecture holds for
sufficiently large $n$. The conjecture also holds for all connected
graphs with diameter $D\leq \frac {2(n-3+2\sqrt{2})}{\pi^2}$ or
minimum degree $\delta\geq \frac n 2$.  We also prove $R\cdot a\geq
\frac {8\sqrt{n-1}}{nD^2}$ and $R\cdot a\geq \frac
{n\delta(2\delta-n+2)} {2(n-1)}$, and then $R\cdot a$ is minimum for
the path if $D\leq
(n-1)^{1/4}$ or $\delta\geq \frac n 2-1$.\\
[2mm] {\bf Keywords:} Randi\'c index; algebraic
connectivity; edge connectivity; diameter; minimum degree\\
[2mm] {\bf AMS Subject Classification 2010:} 05C35, 05C50, 05C90,
92E10.
\end{abstract}

\section{Introduction}
In 1975, Milan Randi\'c \cite{R} proposed a topological index $R$
under the name ``\textit{branching index}", suitable for measuring
the extent of branching of the carbon-atom skeleton of saturated
hydrocarbons. For a graph $G=(V,E)$, the {\it Randi\'c index} $R(G)$
of $G$ was defined as the sum of $1/\sqrt{d(u)d(v)}$ over all edges
$uv$ of $G$, where $d(u)$ denotes the degree of a vertex $u$ in $G$,
i.e., $R(G)=\sum\limits_{uv\in E(G)}\frac 1 {\sqrt{d(u)d(v)}}$\,.
Later, in 1998 Bollob\'as and Erd\"{o}s \cite{BE} generalized this
index by replacing $-\frac{1}{2}$ with any real number $\alpha$,
which is called the general Randi\'c index. Randi\'c  noticed that
there is a good correlation between the Randi\'c index $R$ and
several physico--chemical properties of alkanes: boiling points,
chromatographic retention times, enthalpies of formation, parameters
in the Antoine equation for vapor pressure, surface areas, etc.  For
a comprehensive survey of its mathematical properties, see the book
of Li and Gutman \cite{LG} or the survey \cite{LS}. For terminology
and notations not defined here, we refer the readers to \cite{BM}.

Graffiti is a program designed to make conjectures about, but not
limited to mathematics, in particular graph theory, which was
written by Fajtlowicz from the mid-1980's. A numbered, annotated
listing of several hundred of Graffiti's conjectures can be found in
\cite{F1}. Graffiti has correctly conjectured a number of new bounds
for several well studied graph invariants. The AutoGraphix 1 and 2
systems (AGX 1 and AGX 2) for computer-assisted as well as, for some
functions, fully automated graph theory were developed at GERAD,
Montr\'eal since 1997. The basic idea of the AGX approach is to view
various problems in graph theory, i.e., finding graphs satisfying
some given constraints; finding an extremal graph for some
invariant; finding a conjecture which may be algebraic, i.e., a
relation between graph invariants, or structural, i.e., a
characterization of extremal graphs for some invariant;
Corroborating, refuting and/or strengthening or repairing a
conjecture; and so on. AGX  led to several hundred new conjectures,
ranking from easy ones, proved automatically, to others requiring
longer unassisted or partially assisted proofs, to open ones.

The Laplacian matrix of a graph $G$ is $L(G)=D(G)-A(G)$, where
$D(G)$ is the diagonal matrix of its vertex degrees and $A(G)$ is
the adjacency matrix. Among all eigenvalues of the Laplacian matrix
of $G$, one of the most popular is the second smallest, which was
called the {\it algebraic connectivity} of a graph by Fiedler
\cite{FI} in 1973, denoted by $a(G)$. Its importance is due to the
fact that it is a good parameter to measure, to a certain extent,
how well a graph is connected. For example, it is well-known that a
graph is connected if and only if its algebraic connectivity is
positive. There are several popular graphs for which their algebraic
connectivity is known, such as $a(K_n)=n$, $a(P_n)=2(1-\cos
{\frac{\pi}{n}})$, $a(C_n)=2(1-\cos {\frac{2\pi}{n}})$. Recently,
there is an excellent survey on the algebraic connectivity of graphs
written by de Abreu \cite{Abreu}, for more details, see \cite{C1,
Merris, Mohar, Zhang}.

In \cite{AH}, the authors investigated the relation between Randi\'c
index and the algebraic connectivity of graphs, and proposed the
following two conjectures.

\begin{con}\label{algebraic}
For any connected graph of order $n\geq 3$ with Randi\'c index $R$
and algebraic connectivity $a$, then
\begin{align}\label{eq1} \frac
R a\leq
 \left(\frac{n-3+2\sqrt{2}}{2}\right)/\left(2(1-\cos {\frac{\pi}{n}})\right)
\end{align}
with equality if and only if $G$ is $P_n$.
\end{con}

\begin{con}\label{algebraic1}
For any connected graph on $n\geq 3$ vertices with Randi\'c index
$R$ and algebraic connectivity $a$, then $R\cdot a$ is minimum for

$(1)$ \ path if $n\leq 9$;

$(2)$ \ a balanced double-comet (i.e., two equal sized stars joined
by a path) if $n\geq 10$.
\end{con}

In this paper, we prove that Conjecture \ref{algebraic} holds for
all trees and all connected graphs with edge connectivity
$\kappa'(G)\geq 2$, and if $\kappa'(G)=1$,  Conjecture
\ref{algebraic} holds for sufficiently large $n$. Conjecture
\ref{algebraic} also holds for all connected graphs with diameter
$D\leq \frac {2(n-3+2\sqrt{2})}{\pi^2}$ or $\delta\geq \frac n 2$.
For the lower bound of $R\cdot a$, we show that $R\cdot a\geq \frac
{8\sqrt{n-1}}{nD^2}$ and $R\cdot a\geq \frac {n\delta(2\delta-n+2)}
{2(n-1)}$, and then $R\cdot a$ is minimum for the path if $D\leq
(n-1)^{1/4}$ or $\delta\geq \frac n 2-1$.

\section{Main results}
The following lemma due to \cite{FI} gives a lower bound of the
algebraic connectivity, which is the main tool of our proof.
\begin{lem}\label{alglem1}
For a given connected graph $G$ of order $n$ with the algebraic
connectivity $a$, then $a\geq 2\kappa'(G)(1-\cos{\frac{\pi}{n}})$
and $a\geq 2\delta(G)-n+2$, where $\kappa'(G)$ and $\delta(G)$
denote the edge connectivity and the minimum degree of $G$,
respectively.
\end{lem}

The other lower bound relates to the relationship of the algebraic
connectivity and the diameter, which is from \cite{Mohar1}.
\begin{lem}\label{diameter}
For a graph $G$ of order $n$, its algebraic connectivity $a$ imposes
the upper bound on the diameter $D$ of $G$: $D\geq \frac 4 {na}$.
\end{lem}

The following result comes from \cite{PL}.
\begin{lem}\label{small}
Fix a positive integer $n$. Then among all trees on $n$ vertices the
path has the smallest algebraic connectivity.
\end{lem}

The following bounds of the Randi\'c index is well-known.
\begin{lem}\label{randicindex}
Among all connected graphs of order $n$, regular graphs has the
maximum value of the Randi\'c index, while the star has the minimum
value of the Randi\'c index. Among all trees with $n$ vertices, the
path has the maximum value of the Randi\'c index.
\end{lem}

In \cite{AHZ1}, the authors gave the relation of Randi\'c index and
the minimum degree as follows.

\begin{lem}\label{randicminimumdegree}
For any connected graph on $n\geq 3$ vertices with Randi\'c index
$R$ and minimum degree $\delta$, then  $\frac n {2(n-1)}\leq \frac R
\delta \leq \frac{3n-7+\sqrt 6+3\sqrt 2}{6}$.
\end{lem}

\begin{thm}\label{thm}
For any connected graph of order $n\geq 3$ with Randi\'c index $R$,
algebraic connectivity $a$ and diameter $D$, then

$(1)$ \ the inequality \eqref{eq1} is strict holds for all graphs
with edge connectivity $\kappa'(G)\geq 2$;

$(2)$ \ the inequality \eqref{eq1} holds for all trees, with
equality if and only if the tree is the path;

$(3)$ \ $\frac R a\leq \left(\frac{n}{2}\right)/\left(2(1-\cos
{\frac{\pi}{n}})\right)$ holds for $\kappa'(G)=1$, i.e., the
inequality \eqref{eq1} holds for sufficiently large $n$.
\end{thm}
\pf  $(1)$ \ By Lemma \ref{alglem1}, if the edge connectivity
$\kappa'(G)\geq 2$, then we directly obtain $a\geq
4(1-\cos{\frac{\pi}{n}})$. By Lemma \ref{randicindex}, we know that
for any connected graph $G$, $R\leq \frac n 2$. Therefore, we can
easily to check that for $n\geq 3$,
$$\frac R a<\left(\frac{n}{2}\right)/\left(4(1-\cos
 {\frac{\pi}{n}})\right)< \left(\frac{n-3+2\sqrt{2}}{2}\right)/\left(2(1-\cos
 {\frac{\pi}{n}})\right).$$

Now in the following we assume $\kappa'(G)=1$. Then by Lemma
\ref{alglem1}, we have $a\geq 2(1-\cos{\frac{\pi}{n}})$.

$(2)$ \ Suppose $G$ is a tree. Since path $P_n$ attains the maximum
value of Randi\'c index among all trees, we have $R\leq
\frac{n-3+2\sqrt{2}}{2}$. By Lemma \ref{small}, $a(G)\leq a(P_n)=
2(1-\cos {\frac{\pi}{n}})$. Then
$$\frac R a\leq \left(\frac{n-3+2\sqrt{2}}{2}\right)/\left(2(1-\cos
{\frac{\pi}{n}})\right),$$ with equality if and only if $G$ is
$P_n$.

(3) \ Similar as the proof of (1), we obtain $\frac R a\leq
\left(\frac{n}{2}\right)/\left(2(1-\cos
 {\frac{\pi}{n}})\right)$.

The proof is complete. \qed

\begin{thm}\label{thm}
For any connected graph of order $n\geq 3$ with diameter $D$ and
minimum degree $\delta$, if $D\leq \frac {2(n-3+2\sqrt{2})}{\pi^2}$
or $\delta\geq \frac n 2$, then inequality \eqref{eq1} holds.
\end{thm}
\pf By Lemma \ref{diameter}, we have the lower bound $a\geq \frac 4
{nD}$. If $D\leq \frac {2(n-3+2\sqrt{2})}{\pi^2}$, then we can
easily to check
\begin{align*}
\frac R a&<\left(\frac{n}{2}\right)/\left(\frac{4}{nD}\right)=\frac{n^2D}{8}
\leq \frac{n^2}{8}\cdot \frac{2(n-3+2\sqrt{2})}{\pi^2}\\
&< \left(\frac{n-3+2\sqrt{2}}{2}\right)/\left(2(1-\cos
{\frac{\pi}{n}})\right),
\end{align*}
since $1-\cos {\frac{\pi}{n}}<
1-\left(1-\frac{\pi^2}{n^2}\right)=\frac{\pi^2}{n^2}$.

By Lemmas \ref{alglem1} and \ref{randicminimumdegree}, we know that
$a\geq 2\delta-n+2$ and $R \leq \frac{(3n-7+\sqrt 6+3\sqrt
2)\delta}{6}$. Thus, we have $\frac R a \leq \frac{(3n-7+\sqrt
6+3\sqrt 2)\delta}{6(2\delta-n+2)}$. If $\delta\geq \frac n 2$, then
by some calculations, we can verify that
$$\frac{(3n-7+\sqrt 6+3\sqrt
2)\delta}{6(2\delta-n+2)}<\left(\frac{n-3+2\sqrt{2}}{2}\right)/\left(2(1-\cos
{\frac{\pi}{n}})\right).$$ \qed

\begin{thm}\label{thm}
For any connected graph on $n\geq 3$ vertices with diameter $D$ and
minimum degree $\delta$, then $R\cdot a\geq \frac
{8\sqrt{n-1}}{nD^2}$, $R\cdot a\geq \frac {n\delta(2\delta-n+2)}
{2(n-1)}$, and then $R\cdot a$ is minimum for the path if $D\leq
(n-1)^{1/4}$ or $\delta\geq \frac n 2-1$.

\end{thm}
\pf  By Lemma \ref{randicindex}, the star has the maximum value of
the Randi\'c index among all graphs, then $R\geq \sqrt{n-1}$. Since
$D\geq 2$, we have $R\cdot D\geq 2\sqrt{n-1}$, i.e., $R\geq
\frac{2\sqrt{n-1}}{D}$. Then by Lemma \ref{diameter}, we obtain
$R\cdot a\geq \frac {8\sqrt{n-1}}{nD^2}$. It is easy to check that
if $D\leq (n-1)^{1/4}$, then $\frac {8\sqrt{n-1}}{nD^2}>
(n-3+2\sqrt{2})(1-\cos \frac \pi n)$, i.e., $R\cdot a$ is minimum
for the path if $D\leq (n-1)^{1/4}$.

By Lemmas \ref{alglem1} and \ref{randicminimumdegree}, we know that
$a\geq 2\delta-n+2$ and $R \geq \frac {n\delta} {2(n-1)}$. And then
$R\cdot a\geq \frac {n\delta(2\delta-n+2)} {2(n-1)}$. It is easy to
verify that if $\delta\geq \frac n 2-1$, then $\frac
{n\delta(2\delta-n+2)} {2(n-1)}> (n-3+2\sqrt{2})(1-\cos \frac \pi
n)$, i.e., $R\cdot a$ is minimum for the path if $\delta\geq \frac n
2-1$.  \qed

\end{document}